\documentclass[11pt]{amsart}
\usepackage[T1]{fontenc}
\usepackage[letterpaper,margin=1.12in,headsep=0.24in]{geometry}
\usepackage{amsmath,amssymb,mathtools,microtype,booktabs}
\usepackage{tikz-cd}
\usepackage{mathrsfs,placeins}
\usepackage[hidelinks]{hyperref}
\hypersetup{pdftitle={Fusion rules for group-theoretical fusion 2-categories},pdfsubject={Relative tensor products, stabilizer cochains and the elementary abelian group of order eight},pdfdisplaydoctitle=true}
\newtheorem{theorem}{Theorem}[section]
\newtheorem{proposition}[theorem]{Proposition}
\newtheorem{lemma}[theorem]{Lemma}
\newtheorem{corollary}[theorem]{Corollary}
\theoremstyle{definition}
\newtheorem{definition}[theorem]{Definition}
\theoremstyle{remark}
\newtheorem{remark}[theorem]{Remark}

\newcommand{\A}{\mathbb C^\times}
\newcommand{\dd}{\mathrm d}
\newcommand{\C}{\mathcal C}
\newcommand{\TV}{\mathbf{2Vect}}
\newcommand{\Vect}{\mathbf{Vect}}
\newcommand{\Mod}{\operatorname{Mod}}
\newcommand{\Bimod}{\operatorname{Bimod}}

\newcommand{\CT}{\mathbb{C}^{\times}}

\DeclareMathOperator{\Hom}{Hom}

\title[Fusion rules for group-theoretical fusion 2-categories]{Fusion rules for group-theoretical fusion 2-categories}
\date{}
\usepackage{etoolbox}
\makeatletter
\patchcmd{\@setauthors}{\MakeUppercase{\authors}}{\authors}{}{}
\makeatother
\begin{document}
\raggedbottom

\author[H. Lin]{Hanlin Lin\textsuperscript{a,*}}
\author[Q.-R.~Wang]{Qing-Rui Wang\textsuperscript{b}}
\author[H. Zheng]{Hao Zheng\textsuperscript{c,d}}

\dedicatory{%
  \textsuperscript{a}Qiuzhen College, Tsinghua University,
  Beijing 100084, China\\[3pt]
  \textsuperscript{b}Yau Mathematical Sciences Center, Tsinghua University, Beijing 100084, China\\[3pt]
  \textsuperscript{c}Institute for Applied Mathematics, Tsinghua University,
  Beijing 100084, China\\[3pt]
  \textsuperscript{d}Beijing Institute of Mathematical Sciences and Applications,
  Beijing 101408, China
}

\thanks{E-mail addresses:
  \nolinkurl{linhl24@mails.tsinghua.edu.cn}.}

\begin{abstract}
We give a formula for the relative tensor product of two simple bimodules over $\Vect_H^\psi$ in $\TV_G^\pi$. For each simple summand, the formula specifies a double-coset representative, a stabilizer subgroup, and a two-cochain that trivializes the associator on that subgroup. We use decorated sets to record the bimodule structures and compute this cochain from the coherence isomorphisms of the universal balanced functor. For $G=C_2^3$, we compute the multiplication tables for simple objects in all 36 monoidal equivalence classes. We obtain nine distinct tables and identify the seven Morita equivalence classes. For 25 of the 36 classes, we give explicit models using 2-group-graded 2-vector spaces with a possible four-cocycle twist or 2-representations of finite 2-groups.
\end{abstract}
\maketitle

\section{Introduction}
Let $G$ be a finite group and let $H\le G$ be a subgroup. Let $\pi\in Z^4(G,\CT)$ be a 4-cocycle and let $\psi\in C^3(H,\CT)$ be a 3-cochain satisfying $\dd\psi=\pi|_H$. The category $A=\Vect_H^\psi$ is an algebra in the graded fusion 2-category $\TV_G^\pi$. We write $\C=\C(G,H,\pi,\psi)=\Bimod_{\TV_G^\pi}(A)$ for the 2-category of $A$-$A$ bimodules in $\TV_G^\pi$. A fusion 2-category of this form is called \emph{group-theoretical} \cite{DY}. Our aim is to compute the fusion rules of $\C$ by determining the relative tensor products of simple $A$-$A$ bimodules.

In a $(2+1)$-dimensional quantum field theory, a finite symmetry group $G$ with 't Hooft anomaly $[\pi]$ is implemented by invertible topological surface operators that generate $\TV_G^\pi$. Gauging $H$ with the chosen trivialization $\psi$ of $\pi|_H$ amounts to condensing the algebra $A=\Vect_H^\psi$ in $\TV_G^\pi$. The objects of $\C$ then describe symmetry surfaces in the gauged theory, and bringing two such surfaces together corresponds to taking the relative tensor product of their bimodules over $A$ \cite[\S\S3.3.1 and 3.3.4]{BBFP}. In the untwisted case, Inamura and Ohyama realize these symmetry operators as tensor networks in generalized cluster models \cite{IO}. Our fusion rules determine the simple surfaces in each product and their multiplicities.

D\'ecoppet and Yu identify the underlying finite semisimple 2-category of $\C$ as
\begin{equation}\label{eq:block-decomposition}
  \C\simeq\boxplus_{g\in H\backslash G/H}\Mod(\Vect_{H_g}^{\xi_g}),
  \qquad H_g=H\cap gHg^{-1},
\end{equation}
where $g$ runs over a fixed set of representatives for the double cosets \cite[Cor.~3.2.4]{DY}. By the classification of module categories over $\Vect_{H_g}^{\xi_g}$ \cite{Ostrik}, each simple object of $\C$ can be represented by a triple $(g,L,\phi)$, where $g$ represents a double coset in $H\backslash G/H$, $L\le H_g$ and $\phi$ is a 2-cochain on $L$ satisfying $\dd\phi=\xi_g|_L$. The double-coset ring records the weighted support of fusion \cite[\S3.3.4]{BBFP}, and the component hypergroup records which connected components occur in a product \cite[Cor.~4.16]{ClassF2}. D\'ecoppet and Yu compute the relative tensor product of two free bimodules \cite[Prop.~3.3.1]{DY}. Remark~\ref{rem:known-products} compares the specialization $\pi=1$ with the ordinary relative tensor products of \cite{GP,NikBisets}. A related product occurs in the Burnside ring of finite 2-group representations \cite[Prop.~4.2]{RW}.

Theorem~\ref{thm:triple} gives a formula for the relative tensor product of any two simple objects of $\C$. In particular, it computes the triple labeling each simple summand. The proof derives the six factors contributed by the ambient 4-cocycle $\pi$ and the correction arising from the change of base point in the second input. When both inputs have stabilizer $H$, the formula yields the extension class of the group of invertible objects, as in \eqref{eq:extension-factor}. Our GAP implementation of the fusion formula is available at \url{https://github.com/linlin137/fusion2-gap}.

For $G=C_2^3$, the presentation criterion of \cite[Prop.~3.4.4]{DY} gives 36 monoidal equivalence classes as $H,\pi,\psi$ vary. We compute their multiplication tables in Section~\ref{subsec:cubed} and group the classes by their fusion rules in Table~\ref{tab:cubed-census}. Forgetting $H$ gives seven Morita equivalence classes, listed in Table~\ref{tab:cubed-morita}. Table~\ref{tab:cubed-models} identifies the finite 2-groups for the 19 classes admitting a graded model and the 10 admitting a representation model. Four classes admit both.

\medskip
\noindent\textbf{Acknowledgments.}
This work is supported by the National Natural Science Foundation of China (NSFC) under Grants No. 12274250.

\section{Simple objects and their stabilizer triples}
\label{sec:classification}

Throughout this work, we use normalized multiplicative cochains with values in $\CT$, on which the group acts trivially. For a 2-cochain $\phi$, the inhomogeneous bar differential is given by $(\dd\phi)(a,b,c)=\phi(b,c)\phi(a,bc)/(\phi(ab,c)\phi(a,b))$. We follow the convention of \cite[Def.~2.3.1 and Ex.~2.3.3]{DY}, in which the associator $h_1(h_2h_3)\to(h_1h_2)h_3$ is given by the scalar $\psi(h_1,h_2,h_3)$. 

Fix a set $\mathcal R$ of representatives for $H\backslash G/H$, choosing $1$ for the double coset $H$. A left label $a\in H$ and a right label $b\in H$ preserve degree $g$ when $agb^{-1}=g$. Thus each $a\in H_g=H\cap gHg^{-1}$ determines a unique right label $b=g^{-1}ag$. For three such pairs $(a_i,b_i)$, indexed in the order in which they act, the associator in the $g$-block of \eqref{eq:block-decomposition} has the following value \cite[Prop.~3.2.1, Eq.~(4)]{DY}
\begin{equation}\label{eq:stabilizer-cocycle}
\begin{aligned}
 \xi_g(a_3,a_2,a_1)
 ={}&\frac{\psi(a_3,a_2,a_1)}{\psi(b_1^{-1},b_2^{-1},b_3^{-1})}\\[-1pt]
 &{}\cdot
 \frac{\pi(a_1,a_1^{-1}g,b_2^{-1},b_3^{-1})\,\pi(g,b_1^{-1},b_2^{-1},b_3^{-1})}
      {\pi(a_2,a_1,(a_2a_1)^{-1}g,b_3^{-1})\,
       \pi(a_3,a_2,a_1,(a_3a_2a_1)^{-1}g)} .
\end{aligned}
\end{equation}
\begin{definition}\label{def:triple}
A \emph{stabilizer triple} is a triple $(g,L,\phi)$, where $g\in\mathcal R$, $L\le H_g$ and $\phi\in C^2(L;\A)$ satisfies $\dd\phi=\xi_g|_L$. Write $X(g,L,\phi)$ for the category of right modules over the algebra $\bigoplus_{l \in L}\mathbb{C}_l$ in $\Vect_{H_g}^{\xi_g}$ with multiplication scalars $\phi$. Via \eqref{eq:block-decomposition}, we use the same notation for the corresponding simple $A$-$A$ bimodule.
\end{definition}

Changing the base point in the transitive $H_g$-set conjugates $L$ and changes its multiplication cochain. Write $\operatorname{Ad}_a(u)=aua^{-1}$ for conjugation by $a$. For $a\in H_g$, the correction for this change of base point is
\begin{equation}\label{eq:inner-prism}
 T_a(u,v)=\frac{\xi_g(\operatorname{Ad}_a(u),a,v)}
                  {\xi_g(a,u,v)\,\xi_g(\operatorname{Ad}_a(u),\operatorname{Ad}_a(v),a)}
                  \qquad (u,v\in L).
\end{equation}
It satisfies $\dd T_a=(\operatorname{Ad}_a|_L)^*\xi_g/(\xi_g|_L)$.

\begin{proposition}\label{prop:triple-classification}
Every simple object of $\C(G,H,\pi,\psi)$ is equivalent to an object of the form $X(g,L,\phi)$. Two triples $(g,L,\phi)$ and $(g',L',\phi')$ label equivalent objects if and only if $g=g'$ and there exist $a\in H_g$ and $\lambda\in C^1(L;\A)$ such that
\begin{equation}\label{eq:triple-equivalence}
  L'=aLa^{-1},\qquad \operatorname{Ad}_a^*\phi'=\phi\,T_a\,\dd\lambda .
\end{equation}
For fixed $g$ and $L$, a cochain $\phi$ satisfying $\dd\phi=\xi_g|_L$ exists if and only if $[\xi_g|_L]=[1]$. When it exists, multiplying by two-cocycles gives all choices. Modulo coboundaries these form a torsor over $H^2(L;\A)$, on which $N_{H_g}(L)$ acts by \eqref{eq:triple-equivalence}.
\end{proposition}

\begin{proof}
Apply the classification of indecomposable module categories \cite{Ostrik} to each block in \eqref{eq:block-decomposition}. Changing the base point conjugates the stabilizer and contributes the correction $T_a$ in \eqref{eq:inner-prism}. An isomorphism between the two algebras acts on the homogeneous component of degree $u$ by multiplication by $\lambda(u)$, and compatibility with multiplication gives the factor $\dd\lambda$ in \eqref{eq:triple-equivalence}. These changes account for all equivalences between the labels.
\end{proof}

\section{Computing the product of two triples}\label{sec:triples}

\subsection{The summands}\label{sec:direct-orbits}
Let $S$ and $T$ be the sets of isomorphism classes of simple objects in the graded categories underlying $M=X(g,L,\phi)$ and $N=X(f,L',\phi')$, respectively. For $a,b\in H$ and $s\in S$, write $(a,b)s=asb^{-1}$ for the class obtained by acting on a representative of $s$ on the left by $\mathbb{C}_a$ and on the right by $\mathbb{C}_{b^{-1}}$. We use the same notation for $T$. The two actions commute on isomorphism classes, and inversion in the second coordinate $b$ converts the right action into a left action. Thus this notation defines a left action of $H\times H$ on each set.

The pair $(a,b)$ sends a class of degree $u$ to one of degree $aub^{-1}$. In particular, the pairs preserving degree $g$ are $(a,g^{-1}ag)$ with $a\in H_g$. Under the block description, the degree-$g$ classes in $S$ are identified with the left $H_g$-set $H_g/L$ with base-point stabilizer $L\le H_g$. Similarly, the degree-$f$ classes in $T$ are identified with $H_f/L'$. Choose base points $s_*$ and $t_*$ corresponding to the cosets $1L=L$ and $1L'=L'$, respectively. These cosets have stabilizers $L$ and $L'$ under the respective left actions. The stabilizers of $s_*$ and $t_*$ in $H\times H$ are $\Delta_gL$ and $\Delta_fL'$, respectively, where $\Delta_gL=\{(a,g^{-1}ag):a\in L\}$.

Let $W$ be the set of isomorphism classes of simple objects in the underlying category of $M\boxtimes_A N$. Write $d_S$ and $d_T$ for the degree maps on $S$ and $T$. The middle action $h(s,t)=(sh^{-1},ht)$ of $H$ on $S\times T$ is free, since $sh^{-1}$ has degree $d_S(s)h^{-1}$. Choose skeletal representatives for $M$ and $N$, and write $(s_i,t_i)=h_i(s_{i-1},t_{i-1})$ for successive points in an orbit. In the twisted linearization of its action groupoid, the composite of two basis arrows is the arrow for $h_2h_1$ multiplied by the right and left module reassociation scalars and $\pi(d_S(s_2),h_2,h_1,d_T(t_0))$, by the balancing axiom. For three arrows the two module pentagons cancel the subgroup associators, and $\dd\pi=1$ on $(d_S(s_3),h_3,h_2,h_1,d_T(t_0))$ gives associativity. Every point is then isomorphic to a chosen representative $(s,t)$, whose endomorphism algebra is $\mathbb C$. Evaluation there identifies balanced functors on the orbit with objects of the target in degree $d_S(s)d_T(t)$, and balanced natural transformations with morphisms between those objects. Thus each orbit contributes one copy of $\Vect$ to the relative tensor product \cite[Def.~3.1.3 and Thm.~3.1.6]{Morita}, giving $W=(S\times T)/H$.

The class $[s,t]$ has degree $d_S(s)d_T(t)$, and $(a,c)$ acts on it by $[as,tc^{-1}]$. The orbits of this outer $H\times H$-action index the simple $A$-$A$ bimodule summands of $M\boxtimes_A N$. Equivalently, these summands are indexed by the orbits of the action of $H^3=H\times H\times H$ on $S\times T$ given by $(p,q,r)(s,t)=(psq^{-1},qtr^{-1})$.

Since the $H\times H$-actions on $S$ and $T$ are transitive, every pair in $S\times T$ has the form $(as_*b^{-1},ct_*d^{-1})$. Acting by $(a^{-1},b^{-1},d^{-1})$ sends this pair to $(s_*,b^{-1}ct_*)$. Thus every $H^3$-orbit contains a pair of the form $(s_*,xt_*)$. Two such pairs $(s_*,xt_*)$ and $(s_*,yt_*)$ lie in the same orbit if and only if $y^{-1}bx\in L'$ for some $b\in g^{-1}Lg$. We therefore choose one representative $x\in H$ for each double coset in $(g^{-1}Lg)\backslash H/L'$.

Write $u_x=gxf$ for the degree of $[s_\ast,xt_\ast]$ and $J_x$ for the projection of its outer stabilizer onto the first copy of $H$. Explicitly,
\begin{equation}\label{eq:triple-index}
 x\in(g^{-1}Lg)\backslash H/L',\qquad u_x=gxf,\qquad J_x=L\cap gxL'x^{-1}g^{-1}.
\end{equation}
Indeed, each $a\in J_x$ determines a unique triple $(a,b,c)$ in the stabilizer of $(s_\ast,xt_\ast)$. The right label $b$ at degree $g$, its conjugate $m(a)\in L'$, and the right label $c$ at degree $u_x$ satisfy
\begin{equation}\label{eq:maps}
 b=g^{-1}ag,\qquad m(a)=x^{-1}bx,\qquad c=f^{-1}m(a)f=u_x^{-1}au_x .
\end{equation}
Projection onto the two outer factors identifies this stabilizer with $\Delta_{u_x}J_x$. The cochain $\Phi_x$ below gives the third entry of the label of this summand.

\begin{theorem}\label{thm:triple}
For each representative $x$ chosen above, let $m:J_x\to L'$ be the map in \eqref{eq:maps}. For each degree $v$, let $\theta_v$ be the comparison cochain in \eqref{eq:theta-values}. Let $t_x$ and $\beta_x$ be the cochains in \eqref{eq:thetaexplicit} and \eqref{eq:betaexplicit}, respectively. Set
\begin{equation}\label{eq:Phi}
  \Phi_x=
  \frac{(\phi|_{J_x})\,(m^*\phi')\,(\theta_{u_x}|_{J_x})}
       {(\theta_g|_{J_x})\,(m^*\theta_f)\,t_x\,\beta_x}.
\end{equation}
Then $\dd\Phi_x=\xi_{u_x}|_{J_x}$, and
\begin{equation}\label{eq:mainfusion}
  X(g,L,\phi)\boxtimes_A X(f,L',\phi')
  \simeq\bigoplus_{x\in(g^{-1}Lg)\backslash H/L'}
        X(u_x,J_x,\Phi_x).
\end{equation}
If $u_x\notin\mathcal R$, transport its triple to the chosen representative by \eqref{eq:transport-cochain}.
\end{theorem}

For a fixed $x$ we often write $u=u_x$ and $J=J_x$. The proof of Theorem~\ref{thm:triple} is completed in Section~\ref{sec:product-proof}.
\subsection{The comparison cochain of a simple bimodule}\label{sec:comparison-scalars}
To compute the third entry of a product triple, we use a second description of each input. The triple $(g,L,\phi)$ records the multiplication at the chosen point $s_\ast$. Equivalently, we can keep the whole transitive graded $H$-biset $S$ and record how comparisons compose at every point of $S$. We denote the resulting two-cochain by $b_S$ and explain below how to pass between the two descriptions.

For an $H$-biset $S$, write $S//(H\times H)$ for its action groupoid. Its objects are the elements of $S$, and an arrow labeled by $(a,b)$ goes from $s$ to $asb^{-1}$. A $G$-grading on $S$ is a map $d_S:S\to G$ satisfying $d_S(asb^{-1})=a\,d_S(s)b^{-1}$. It induces a functor $S//(H\times H)\to G//(H\times H)$, also denoted by $d_S$.

\begin{definition}\label{def:decorated-set}
A \emph{decorated set} for $(G,H,\pi,\psi)$ is a finite $G$-graded $H$-biset $S$ together with a normalized two-cochain $b_S\in C^2(S//(H\times H),\CT)$ satisfying
\begin{equation}\label{eq:decoration}
 \dd b_S=d_S^*\alpha .
\end{equation}
To define $\alpha\in Z^3(G//(H\times H);\A)$, take any $v_0\in G$ and three successive arrows $k_i=(a_i,b_i)$. Set $v_i=a_i v_{i-1}b_i^{-1}$ for $i=1,2,3$. Then
\begin{equation}\label{eq:alphaexpanded}
 \alpha(v_0;k_3,k_2,k_1)=
 \frac{\pi(a_3,a_2,a_1,v_0)}{\pi(a_3,a_2,v_1,b_1)}
 \frac{\pi(a_3,v_2,b_2,b_1)}{\pi(v_3,b_3,b_2,b_1)}
 \frac{\psi(b_3,b_2,b_1)}{\psi(a_3,a_2,a_1)} .
\end{equation}

For four successive arrows, the four-cocycle identity for $\pi$ and the relations $a_i v_{i-1}=v_i b_i$ give $\pi(a_4,a_3,a_2,a_1)/\pi(b_4,b_3,b_2,b_1)$ as the differential of the product of the $\pi$-factors in \eqref{eq:alphaexpanded}. The differential of the $\psi$-ratio is its inverse, since $\dd\psi=\pi|_H$. Thus $\dd\alpha=1$.

Two decorated sets are equivalent if a degree-preserving biset isomorphism identifies their cochains up to a coboundary.
\end{definition}

\begin{remark}\label{rem:biset-decoration}
When $\pi=\psi=1$, the underlying $(H,H)$-biset is decorated by the class of $b_S$ in $H^2(H\times H;\operatorname{Fun}(S,\A))$ for the permutation action \cite[Def.~2.1]{NikBisets}. Equivalently, this decoration is specified by the restriction of $b_S$ to one point stabilizer in each orbit, up to coboundary \cite[Lem.~3.1]{GRR}.
\end{remark}

For an $A$-$A$ bimodule $M$, let $S$ be the set of isomorphism classes of simple objects in its underlying category. Take any $s_0\in S$ and two successive arrows $k_i=(a_i,b_i)$, so that $s_i=a_i s_{i-1}b_i^{-1}$ for $i=1,2$. Choose a representative of every simple class and an isomorphism $q_i:a_i s_{i-1}\to s_i b_i$ for each arrow, using the identity for an identity arrow. Write $q_{21}$ for the chosen comparison of $k_2k_1$.

Write $\ell^M_{a_2,a_1;s_0}:(a_2a_1)s_0\to a_2(a_1s_0)$, $\mu^M_{a_2;s_1;b_1}:(a_2s_1)b_1\to a_2(s_1b_1)$, and $r^M_{s_2;b_2,b_1}:(s_2b_2)b_1\to s_2(b_2b_1)$ for the left, mixed, and right associators of $M$, respectively. We can compare $q_{21}$ with the composite
$$
\begin{aligned}
 (a_2a_1)s_0 &\stackrel{\ell}{\longrightarrow} a_2(a_1s_0)
 \stackrel{q_1}{\longrightarrow} a_2(s_1b_1)\stackrel{\mu^{-1}}{\longrightarrow} (a_2s_1)b_1\\
 &\stackrel{q_2}{\longrightarrow} (s_2b_2)b_1\stackrel{r}{\longrightarrow} s_2(b_2b_1).
\end{aligned}
$$
The source and target are isomorphic simple objects, so the direct comparison $q_{21}$ is a unique scalar multiple of this composite. In the scalar formulas below we work with skeletal representatives of the simple objects and suppress the identity maps on unchanged factors.

\begin{lemma}\label{lem:comparison-scalar}
The comparison scalars
\begin{equation}\label{eq:input-comparison}
 b_S(s_0;k_2,k_1)=
 \frac{q_{21}\,\mu^M_{a_2;s_1;b_1}}
 {q_2 q_1\,\ell^M_{a_2,a_1;s_0}\,r^M_{s_2;b_2,b_1}} 
\end{equation}
determine a decoration on $S$. Replacing $q$ by $q\lambda$ replaces $b_S$ by $b_S/\dd\lambda$.
\end{lemma}

\begin{proof}
For three successive arrows in $S//(H\times H)$, we must show that
\begin{equation}\label{eq:direct-composition}
 \frac{b_S(s_0;k_2,k_1)\,b_S(s_0;k_3,k_2k_1)}
      {b_S(s_0;k_3k_2,k_1)\,b_S(s_1;k_3,k_2)}
 =\alpha(v_0;k_3,k_2,k_1).
\end{equation}
Write $C(u,v,w):(uv)w\to u(vw)$ for forward reassociation. With one module entry $C$ is $\ell$, $\mu$ or $r$ according to its position. On three subgroup entries it has scalar $\psi^{-1}$. For homogeneous entries, with $|u|$ denoting the degree of $u$, the pentagon gives
\begin{equation}\label{eq:pent}
 E(u,v,w,x)=\frac{C(u,v,w)\,C(u,vw,x)\,C(v,w,x)}
 {C(uv,w,x)\,C(u,v,wx)}=\pi(|u|,|v|,|w|,|x|)^{-1}.
\end{equation}
Let $E_0,E_1,E_2,E_3$ be its values on
$$
 (a_3,a_2,a_1,s_0),\quad(a_3,a_2,s_1,b_1),\quad
 (a_3,s_2,b_2,b_1),\quad(s_3,b_3,b_2,b_1).
$$
In $E_1E_3/(E_0E_2)$ the common faces $\ell_{a_3,a_2;a_1s_0}$, $\mu_{a_3;a_2s_1;b_1}$ and $r_{a_3s_2;b_2,b_1}$ cancel because $a_i s_{i-1}=s_i b_i$ on simple classes. Substituting \eqref{eq:input-comparison} into the left side of \eqref{eq:direct-composition} cancels all the $q$-factors and gives this ratio times $\psi(b_3,b_2,b_1)/\psi(a_3,a_2,a_1)$. This is \eqref{eq:alphaexpanded}.
\end{proof}

\begin{lemma}\label{lem:base-point-transport}
Let $b$ be a decoration satisfying \eqref{eq:decoration}, and let $e:s\to s'$ be an arrow in the action groupoid. Write $g=d_S(s)$ and let $\Delta_gL$ be the stabilizer of $s$, with $L\le H_g$. Identify $L$ with $\Delta_gL$ via $a\mapsto A_a=(a,g^{-1}ag)$, and put $A'_a=eA_ae^{-1}$. For $a,a'\in L$, define
\begin{equation}\label{eq:prism-transport}
  T_e(a',a)=
  \frac{\alpha(g;A'_{a'},e,A_a)}
       {\alpha(g;A'_{a'},A'_a,e)\,\alpha(g;e,A_{a'},A_a)}.
\end{equation}
Set $\kappa_e(a)=b(s;A'_a,e)/b(s;e,A_a)$. Then
\begin{equation}\label{eq:direct-transport}
  b(s';A'_{a'},A'_a)
  =b(s;A_{a'},A_a)\,T_e(a',a)\,(\dd\kappa_e)(a',a).
\end{equation}
\end{lemma}

\begin{proof}
Expand the three factors in \eqref{eq:prism-transport} using \eqref{eq:direct-composition}. The identities $eA_a=A'_ae$ and $eA_{a'}=A'_{a'}e$ cancel the terms containing these composites. Using $A_{a'}A_a=A_{a'a}$ and $A'_{a'}A'_a=A'_{a'a}$, the remaining terms give the ratio $b(s';A'_{a'},A'_a)/b(s;A_{a'},A_a)$ multiplied by $\kappa_e(a'a)/(\kappa_e(a')\kappa_e(a))$. The latter factor is $(\dd\kappa_e)(a',a)^{-1}$, which gives \eqref{eq:direct-transport}.
\end{proof}

Let $\ell_g:BH_g\to G//(H\times H)$ send $a\in H_g$ to the loop $(a,g^{-1}ag)$ at $g$. The two-cochain $\theta_g$ of Appendix~\ref{app:comparison} satisfies
\begin{equation}\label{eq:xi}
 \xi_g=\frac{\dd\theta_g}{\ell_g^*\alpha}.
\end{equation}

\begin{theorem}\label{thm:decorated-classification}
Transitive decorated sets up to equivalence correspond bijectively to stabilizer triples up to \eqref{eq:triple-equivalence}. For $s_*$ of degree $g\in\mathcal R$ with stabilizer $\Delta_g L$, the corresponding triple is $(g,L,\phi)$ with
\begin{equation}\label{eq:dictionary}
 \phi\equiv\frac{\theta_g|_L}{b_S|_{s_*,L}}\pmod{B^2(L;\A)} .
\end{equation}
\end{theorem}

\begin{proof}
Let $(S,b_S)$ be a transitive decorated set. Choose $s_*\in S$ of degree $g\in\mathcal R$ with stabilizer $\Delta_g L$. Then $S\simeq(H\times H)/\Delta_g L$ as a graded biset, and the image of $d_S$ is $HgH$. Restricting $\alpha$ to this component, we have a commutative diagram
$$
\begin{tikzcd}
BL \arrow[d, "\iota"'] \arrow[r, "i"]
  & S//(H\times H) \arrow[d, "d_S"] \\
BH_g \arrow[r, "\ell_g"']
  & (HgH)//(H\times H).
\end{tikzcd}
$$
Here $\iota$ is induced by the inclusion $L\le H_g$. Both horizontal functors are equivalences, since they include the full subgroupoid on a single object in a connected groupoid. Put $b_L=i^*b_S$ and define $\phi=\theta_g|_L/b_L$. The commutative square gives $i^*d_S^*\alpha=\iota^*\ell_g^*\alpha=(\ell_g^*\alpha)|_L$. Hence \eqref{eq:xi} gives
$$
\dd\phi
=\frac{\xi_g|_L\,(\ell_g^*\alpha)|_L}{i^*d_S^*\alpha}
=\xi_g|_L.
$$
Multiplying $b_S$ by a coboundary changes $\phi$ by a coboundary, so the assignment is well defined on classes.

Conversely, if $\dd\phi=\xi_g|_L$, then $b_L=\theta_g|_L/\phi$ satisfies $\dd b_L=i^*d_S^*\alpha$. Since $i^*$ is an isomorphism in degree three and $i^*d_S^*\alpha$ is a coboundary, $d_S^*\alpha$ is a coboundary as well. Thus $\dd b_S=d_S^*\alpha$ has a solution. Any two solutions differ by a two-cocycle. Since $i^*$ is also an isomorphism in degree two, exactly one class of solutions restricts to $b_L$ modulo coboundaries. This proves the bijection for the chosen base point.

Changing the base point by $e=(p,g^{-1}pg)$ replaces $L$ by $L'=\operatorname{Ad}_p(L)$. By \eqref{eq:direct-transport}, $\operatorname{Ad}_p^*b_{L'}\equiv b_LT_e$ modulo coboundaries. The identity following \eqref{eq:transport-cochain} gives $\operatorname{Ad}_p^*(\theta_g|_{L'})\equiv(\theta_g|_L)T_eT_p$. Taking the quotient gives $\operatorname{Ad}_p^*\phi'\equiv\phi T_p$, as in \eqref{eq:triple-equivalence}. Thus the bijection descends to equivalence classes.

To identify the multiplication at $s_*$, take two stabilizer pairs $(a_i,b_i)$ with $a_i\in L$ and $b_i=g^{-1}a_i g$. In Lemma~\ref{lem:comparison-scalar}, the composite built from $q_1:a_1s_*\to s_*b_1$ and $q_2:a_2s_*\to s_*b_2$ has coefficient $b_L(a_2,a_1)^{-1}$ relative to $q_{21}$. In the stabilizer convention, $q_i$ gives a map $(a_i s_*)b_i^{-1}\to s_*$ by cancelling $b_i b_i^{-1}$. Appendix~\ref{app:comparison} compares cancelling the two right factors successively with cancelling their product. With the choice \eqref{eq:theta-values}, the multiplication coefficient in $\Vect_{H_g}^{\xi_g}$ is $\theta_g(a_2,a_1)b_L(a_2,a_1)^{-1}$. Rescaling the $q_i$ changes this coefficient by a coboundary, as in \eqref{eq:dictionary}.
\end{proof}

\begin{remark}
Changing $\theta_g$ by a two-cocycle $z_g$ multiplies $\phi$ by $z_g|_L$. In \eqref{eq:Phi} the changes cancel for the two inputs and multiply the output cochain by $z_u|_J$.
\end{remark}

\subsection{The comparison scalar on the balanced product}

Let $S$, $T$, and $W$ denote the sets of isomorphism classes of simple objects in the underlying categories of $M=X(g,L,\phi)$, $N=X(f,L',\phi')$, and $M\boxtimes_A N$, respectively. We compute $b_W$ from $b_S$ and $b_T$. Let $F$ be the universal balanced functor $M\times N \to M\boxtimes_A N$, and consider two successive arrows $\gamma_i=(a_i,c_i):w_{i-1}\to w_i$ in $W//(H\times H)$. Choose one representative for each orbit in $W$, and write $(s_i,t_i)$ for the representative of $w_i$. Since the middle $H$-action on $S\times T$ is free, there is a unique $h_i\in H$ such that
\begin{equation}
 a_i s_{i-1}=s_i h_i,\qquad h_i t_{i-1}=t_i c_i.
 \label{eq:labels}
\end{equation}
Write $F_i=F(s_i,t_i)$ and $x_{21}=x_2x_1$ for $x=a,h,c$. Applying \eqref{eq:labels} twice shows that $h_{21}$ is the middle label for $\gamma_2\gamma_1$. Write $q_i^S:a_i s_{i-1}\to s_i h_i$ and $q_i^T:h_i t_{i-1}\to t_i c_i$ for the input comparisons and $q_{21}^S,q_{21}^T$ for those of the composite. The three structure isomorphisms of $F$ have the form
\begin{align*}
 L(a,s,t)&: F(as,t)\longrightarrow aF(s,t),\\
 \kappa(s,h,t)&: F(sh,t)\longrightarrow F(s,ht),\\
 \rho(s,t,c)&: F(s,t)c\longrightarrow F(s,tc).
\end{align*}
Starting at $a_i F_{i-1}$, we apply $L_i^{-1}$ and $F(q_i^S,1)$ to reach $F(s_i h_i,t_{i-1})$. Balancing moves $h_i$ to the second input. We then apply $F(1,q_i^T)$ and $\rho_i^{-1}$ to reach $F_i c_i$. Here $L_i=L(a_i,s_{i-1},t_{i-1})$, $\kappa_i=\kappa(s_i,h_i,t_{i-1})$ and $\rho_i=\rho(s_i,t_i,c_i)$. Using $s_0,t_0$ and $s_2,t_2$ as the initial and final representatives for the composite, we obtain
\begin{equation}
 V_i=\frac{q_i^S\,\kappa_i\,q_i^T}{L_i\,\rho_i},\qquad i=1,2,21.
 \label{eq:V}
\end{equation}
In the following table, the subscript $W$ on a constraint refers to $M\boxtimes_A N$. The columns list the left, mixed, and right constraints.
\[
\begin{array}{c@{\qquad}ccc}
\toprule
 X & \ell_X & \mu_X & r_X\\ \midrule
 M&\ell^M_{a_2,a_1;s_0}&\mu^M_{a_2;s_1;h_1}&r^M_{s_2;h_2,h_1}\\[2pt]
 N&\ell^N_{h_2,h_1;t_0}&\mu^N_{h_2;t_1;c_1}&r^N_{t_2;c_2,c_1}\\[2pt]
 W&\ell^W_{a_2,a_1;F_0}&\mu^W_{a_2;F_1;c_1}&r^W_{F_2;c_2,c_1}\\
\bottomrule
\end{array}
\]
As in \eqref{eq:input-comparison}, the two-step comparison applies $\ell_W$, $a_2V_1$, $\mu_W^{-1}$, $V_2c_1$ and $r_W$ in that order. Substitution of \eqref{eq:V} in the quotient of the direct and two-step comparisons gives
\begin{equation}
\begin{aligned}
 b_W={}&\frac{q^S_{21}\,\kappa_{21}\,q^T_{21}}{L_{21}\,\rho_{21}}
         \frac{L_1\,\rho_1}{q^S_1\,\kappa_1\,q^T_1}
         \frac{L_2\,\rho_2}{q^S_2\,\kappa_2\,q^T_2}
         \frac{\mu_W}{\ell_W\,r_W}.
\end{aligned}
\label{eq:raw}
\end{equation}
The same definition applied separately to the two inputs gives
\begin{align}
 b_S&=\frac{q^S_{21}\,\mu_M}{q^S_1 q^S_2\,\ell_M r_M},\label{eq:bs}\\
 b_T&=\frac{q^T_{21}\,\mu_N}{q^T_1 q^T_2\,\ell_N r_N}.\label{eq:bt}
\end{align}
Dividing \eqref{eq:raw} by these two expressions cancels all six input comparisons. If we abbreviate $U_1U_2/U_{21}$ by $\delta U$, the remaining factors are
\begin{equation}
 \mathcal E:=\frac{b_W}{b_S b_T}
 =\frac{\delta L\,\delta\rho}{\delta\kappa}
  \frac{\mu_W}{\ell_W r_W}
  \frac{\ell_M r_M}{\mu_M}
  \frac{\ell_N r_N}{\mu_N}.
\label{eq:res}
\end{equation}
To evaluate this quotient, we use the six coherence diagrams for $F$ \cite{Morita}. Write $st=F(s,t)$ in the reassociation $C$ of \eqref{eq:pent}, so that $C(a,s,t)=L(a,s,t)$, $C(s,h,t)=\kappa(s,h,t)$ and $C(s,t,c)=\rho(s,t,c)$. The left and right module-functor axioms give (P1) and (P6), and their mixed compatibility gives (P4). Compatibility of balancing with the exterior actions gives (P2) and (P5). Balancing two successive middle actions gives (P3). With $g_i=d_S(s_i)$ and $f_i=d_T(t_i)$, these equations read
\begin{align}
 P_1 &=\pi(a_2,a_1,g_0,f_0)^{-1}
 =\frac{\ell_M\,L(a_2,a_1s_0,t_0)\,L_1}{L_{21}\,\ell_W},
 \tag{P1}\label{p1}\\[3pt]
 P_2 &=\pi(a_2,g_1,h_1,f_0)^{-1}
 =\frac{\mu_M\,L(a_2,s_1h_1,t_0)\,\kappa_1}
       {\kappa(a_2s_1,h_1,t_0)\,L(a_2,s_1,h_1t_0)},
 \tag{P2}\label{p2}\\[3pt]
 P_3 &=\pi(g_2,h_2,h_1,f_0)^{-1}
 =\frac{r_M\,\kappa_{21}\,\ell_N}
       {\kappa(s_2h_2,h_1,t_0)\,\kappa(s_2,h_2,h_1t_0)},
 \tag{P3}\label{p3}\\[3pt]
 P_4 &=\pi(a_2,g_1,f_1,c_1)^{-1}
 =\frac{L_2\,\mu_W\,\rho_1}{\rho(a_2s_1,t_1,c_1)\,L(a_2,s_1,t_1c_1)},
 \tag{P4}\label{p4}\\[3pt]
 P_5 &=\pi(g_2,h_2,f_1,c_1)^{-1}
 =\frac{\kappa_2\,\rho(s_2,h_2t_1,c_1)\,\mu_N}
       {\rho(s_2h_2,t_1,c_1)\,\kappa(s_2,h_2,t_1c_1)},
 \tag{P5}\label{p5}\\[3pt]
 P_6 &=\pi(g_2,f_2,c_2,c_1)^{-1}
 =\frac{\rho_2\,\rho(s_2,t_2c_2,c_1)\,r_N}{r_W\,\rho_{21}} .
 \tag{P6}\label{p6}
\end{align}
Each pentagon has two module entries. This explains why none of the six equations has a factor from the subgroup associator. To recover \eqref{eq:res}, multiply (P1), (P3), (P4) and (P6) and divide by (P2) and (P5). For instance, $a_1s_0=s_1h_1$ makes the $L$-factors from the first two rows cancel. The six cancellations all follow from \eqref{eq:labels}.
\[
\begin{array}{@{}cll@{}}
\toprule
 \text{Rows}& \text{Numerator factor}&\text{Denominator factor}\\ \midrule
 1,2 & L(a_2,a_1s_0,t_0)        & L(a_2,s_1h_1,t_0)\\
 2,4 & L(a_2,s_1,h_1t_0)        & L(a_2,s_1,t_1c_1)\\
 2,3 & \kappa(a_2s_1,h_1,t_0)   & \kappa(s_2h_2,h_1,t_0)\\
 5,3 & \kappa(s_2,h_2,t_1c_1)   & \kappa(s_2,h_2,h_1t_0)\\
 5,4 & \rho(s_2h_2,t_1,c_1)     & \rho(a_2s_1,t_1,c_1)\\
 6,5 & \rho(s_2,t_2c_2,c_1)     & \rho(s_2,h_2t_1,c_1)\\
\bottomrule
\end{array}
\]
Keeping the remaining factors in row order gives
\begin{align*}
 \frac{P_1 P_3 P_4 P_6}{P_2 P_5}
  ={}&\underbrace{\frac{\ell_M L_1}{L_{21}\ell_W}}_{\text{from }P_1}
      \underbrace{\frac{1}{\mu_M\kappa_1}}_{\text{from }P_2^{-1}}
      \underbrace{(r_M\kappa_{21}\ell_N)}_{\text{from }P_3}\\
   &{}\cdot\underbrace{(L_2\mu_W\rho_1)}_{\text{from }P_4}
      \underbrace{\frac{1}{\kappa_2\mu_N}}_{\text{from }P_5^{-1}}
      \underbrace{\frac{\rho_2 r_N}{r_W\rho_{21}}}_{\text{from }P_6}\\
  ={}&\frac{\delta L\,\delta\rho}{\delta\kappa}
      \frac{\mu_W}{\ell_W r_W}
      \frac{\ell_M r_M}{\mu_M}\frac{\ell_N r_N}{\mu_N}
      =\mathcal E .
\end{align*}
Substituting the six pentagonators into $b_W=b_S b_T\mathcal E$ gives
\begin{equation}
\begin{aligned}
 b_W(w_0;\gamma_2,\gamma_1)
  ={}& b_S(s_0;(a_2,h_2),(a_1,h_1))\\
    &{}\cdot b_T(t_0;(h_2,c_2),(h_1,c_1))\\
    &{}\cdot\frac{\pi(a_2,g_1,h_1,f_0)}
      {\pi(a_2,a_1,g_0,f_0)\,\pi(g_2,h_2,h_1,f_0)}\\
    &{}\cdot\frac{\pi(g_2,h_2,f_1,c_1)}
      {\pi(a_2,g_1,f_1,c_1)\,\pi(g_2,f_2,c_2,c_1)} .
\end{aligned}
\label{eq:direct-product}
\end{equation}
The calculation used only representatives of the three classes along these arrows. It therefore also applies directly to pairs in $S\times T$. As in \eqref{eq:input-comparison}, rescaling the output comparisons changes the decoration by a two-coboundary.

\subsection{Proof of Theorem \ref{thm:triple}}\label{sec:product-proof}

The summand indexed by $x$ is represented by $(s_\ast,xt_\ast)$, whereas $\phi'$ describes the second input at $t_\ast$. We apply Lemma~\ref{lem:base-point-transport} to the second input along $(x,1):t_\ast\to xt_\ast$. For $a_1,a_2\in J$, let $b_i$ and $c_i$ be the coordinates determined by $a_i$ in \eqref{eq:maps}. The pair $(m(a_i),c_i)$ fixes $t_\ast$, and $xm(a_i)=b_i x$ gives $(x,1)(m(a_i),c_i)=(b_i,c_i)(x,1)$. Thus $(b_i,c_i)$ fixes $xt_\ast$. Since $t_\ast$ has degree $f$, the correction $t_x=T_{(x,1)}$ is given by
\begin{equation}\label{eq:thetaexplicit}
  t_x(a_2,a_1)=
  \frac{\alpha(f;(b_2,c_2),(x,1),(m(a_1),c_1))}
       {\alpha(f;(b_2,c_2),(b_1,c_1),(x,1))}
  \frac{1}{\alpha(f;(x,1),(m(a_2),c_2),(m(a_1),c_1))}.
\end{equation}
Equation~\eqref{eq:direct-transport} therefore gives the congruence
\begin{equation*}
\begin{aligned}
  &b_T(xt_\ast;(b_2,c_2),(b_1,c_1))\\
  &\qquad\equiv
  b_T(t_\ast;(m(a_2),c_2),(m(a_1),c_1))\,t_x(a_2,a_1)
  \pmod{B^2(J;\A)}.
\end{aligned}
\end{equation*}

The degrees of $s_\ast$ and $xt_\ast$ are $g$ and $xf$, respectively. Evaluating the six $\pi$-factors in \eqref{eq:direct-product} on the stabilizer of $(s_\ast,xt_\ast)$ gives
\begin{equation}\label{eq:betaexplicit}
\begin{aligned}
  \beta_x(a_2,a_1)={}&
  \frac{\pi(a_2,g,b_1,xf)}
       {\pi(a_2,a_1,g,xf)\,\pi(g,b_2,b_1,xf)}
  \frac{\pi(g,b_2,xf,c_1)}
       {\pi(a_2,g,xf,c_1)\,\pi(g,xf,c_2,c_1)}.
\end{aligned}
\end{equation}
\begin{proof}[Proof of Theorem~\ref{thm:triple}]
By \eqref{eq:dictionary}, the first input contributes $(\theta_g|_J)/(\phi|_J)$ to the decoration at $(s_\ast,xt_\ast)$. For the second input, pull $(\theta_f|_{L'})/\phi'$ back along $m:J\to L'$ and use \eqref{eq:direct-transport} to move from $t_\ast$ to $xt_\ast$. This gives $(m^*\theta_f)t_x/(m^*\phi')$ up to a coboundary. By \eqref{eq:direct-product}, the output decoration on $J$ is the product of these two cochains and $\beta_x$. Dividing $\theta_u|_J$ by this decoration gives exactly $\Phi_x$ in \eqref{eq:Phi} as a representative of the output label. Its differential is $(\dd\theta_u)|_J/((\ell_u^*\alpha)|_J)=\xi_u|_J$ by \eqref{eq:direct-composition} and \eqref{eq:xi}. Taking the direct sum over the orbits in \eqref{eq:triple-index} proves \eqref{eq:mainfusion}.
\end{proof}

We now express \eqref{eq:mainfusion} at the chosen representative of the output double coset. If $v=puq^{-1}\in \mathcal{R}$ is a chosen double coset representative, the arrow $e=(p,q)\in H\times H$ takes the loops $A_a=(a,u^{-1}au)$ at degree $u$ to loops at degree $v$. Their first labels give $j:J\to J'$, where $j=\operatorname{Ad}_p|_J$ and $J'=j(J)$. Applying \eqref{eq:direct-transport} to the output decoration and then using \eqref{eq:dictionary} gives
\begin{equation}\label{eq:transport-cochain}
 j^*\Phi_v\equiv
 \frac{\Phi_x\,j^*(\theta_v|_{J'})}{T_e\,(\theta_u|_J)}
 \pmod{B^2(J;\A)} .
\end{equation}
This transport also accounts for the choices made in indexing the summands. If $x'=kx\ell$ with $k\in g^{-1}Lg$ and $\ell\in L'$, the arrow $(gkg^{-1},k,f^{-1}\ell^{-1}f)$ sends $(s_\ast,xt_\ast)$ to $(s_\ast,x't_\ast)$. Thus \eqref{eq:transport-cochain} identifies the labels from these two representatives. To compare two choices of transport to a given degree, it suffices to consider a loop at the source. For such a loop $e=(p,u^{-1}pu)$, put $I_p\theta_u(a)=\theta_u(p,a)/\theta_u(j(a),p)$. Inserting $p$ successively in the three positions of $\dd\theta_u=\xi_u\,\ell_u^*\alpha$ gives $j^*\theta_u/(\theta_u T_e)=T_p/\dd(I_p\theta_u)$, with $T_p$ as in \eqref{eq:inner-prism}. The difference between the two labels is therefore precisely the equivalence of Proposition~\ref{prop:triple-classification}.

Now choose input base points $s'_\ast=a s_\ast b^{-1}$ and $t'_\ast=c t_\ast d^{-1}$. Setting $x'=bxc^{-1}$, we see that $(a,b,d)$ sends $(s_\ast,xt_\ast)$ to $(s'_\ast,x't'_\ast)$. This identifies the two double-coset index sets, and the exterior labels $(a,d)$ transport the output by \eqref{eq:transport-cochain}. Multiplying the input cochains by $\dd\eta$ and $\dd\eta'$ multiplies the output cochain by $\dd((\eta|_J)m^*\eta')$. Consequently, the multiplicity of a fixed simple object is the number of double cosets that give its label up to equivalence. In the following formula, each output triple is first transported to its chosen representative by \eqref{eq:transport-cochain}.
\begin{equation}\label{eq:multiplicity}
 N_{(g,L,\phi),(f,L',\phi')}^{(k,K,\varphi)}
 =\sum_{x\in(g^{-1}Lg)\backslash H/L'}
        \mathbf{1}\{(u_x,J_x,\Phi_x)\sim(k,K,\varphi)\}.
\end{equation}
\begin{remark}\label{rem:known-products}
For $\pi=1$, Galindo and Plavnik describe the ungraded product over $\Vect_H^\psi$ by equivariantization \cite[Thm.~7.15 and Rem.~7.17(2)]{GP}. If also $\psi=1$, then $\Phi_x=(\phi|_{J_x})\,m^*\phi'$. By \eqref{eq:dictionary}, $[\Phi_x]^{-1}$ is the output decoration in Nikshych's double-coset formula for trivial middle stabilizer \cite[Thm.~2.5]{NikBisets}. Taking $H=G$ gives $J_x=L\cap xL'x^{-1}$ and recovers the fusion rule for simple $2$-representations of $G$ \cite[\S4.3.3]{BBFPI}. For general $\pi$, \eqref{eq:Phi} includes the ambient phase $\beta_x$ and the base-point correction $t_x$ before transport to the chosen output degree.
\end{remark}

\section{Applications of the fusion rule}\label{sec:examples}

\subsection{The group of invertible objects}
Let $K_0^\oplus(\C)_+$ denote the split Grothendieck semiring of $\C$. Its elements are the equivalence classes of objects, with addition induced by direct sum and multiplication induced by $\boxtimes_A$. Theorem~\ref{thm:triple} gives both its multiplication and a way to recognize which basis elements are invertible. We begin by identifying the unit. The restriction $\ell_1^*\alpha$ to diagonal loops at $1$ is trivial, so $\dd\theta_1=\xi_1$. Thus $U=X(1,H,\theta_1)$ is defined. For any simple object $X$, the products $UX$ and $XU$ each have a single summand, indexed by $x=1$. Normalization gives $t_1=\beta_1=1$ in these products, and \eqref{eq:Phi} gives back the cochain of $X$. Hence $U$ is the unit.

A \emph{positive character} of $K_0^\oplus(\C)_+$ is a unital semiring homomorphism $\delta:K_0^\oplus(\C)_+\to\mathbb R_{\geq 0}$ such that $\delta([X])>0$ for every simple object $X$ of $\C$. We write $\delta(X)$ for $\delta([X])$.

\begin{corollary}\label{cor:rank}
The semiring $K_0^\oplus(\C)_+$ has a unique positive character $d$, given on simple objects by $d(X(g,L,\phi))=|H|/|L|$.
\end{corollary}
\begin{proof}
Let $X=X(g,L,\phi)$ and $X'=X(f,L',\phi')$. The double coset $(g^{-1}Lg)xL'$ has $|L||L'|/|J_x|$ elements. Summing over the double cosets in $H$ gives $\sum_x |H|/|J_x|=(|H|/|L|)(|H|/|L'|)$, as required for multiplicativity. The unit has stabilizer $H$, so $d(U)=1$.

To prove uniqueness, let $\bar f\in\mathcal R$ be the chosen representative of $HfH$ for each $f\in G$, and set $R=\bigoplus_{f\in G}X(\bar f,1,1)$. Fix $X=X(g,L,\phi)$ and choose representatives $x\in H$ for the cosets in $(g^{-1}Lg)\backslash H$. Taking a base point of degree $f$ in $X(\bar f,1,1)$, the fusion formula gives one summand of degree $gxf$ with trivial stabilizer for each $x$. For fixed $x$, the map $f\mapsto gxf$ permutes $G$, so summing over $f$ contributes one copy of $R$. There are $|H|/|L|$ choices of $x$, hence $XR=(|H|/|L|)R$ in the split Grothendieck semiring. If $\delta$ is a positive character, then $\delta(R)>0$. Applying $\delta$ to this identity and dividing by $\delta(R)$ gives $\delta(X)=|H|/|L|$.
\end{proof}
\begin{lemma}\label{lem:invertible-triples}
A simple object $X(g,L,\phi)$ is invertible if and only if $L=H$.
\end{lemma}

\begin{proof}
If $X=X(g,L,\phi)$ is invertible, then $X^{-1}$ is also simple. The positive integers $d(X)$ and $d(X^{-1})$ have product $1$. Thus $d(X)=1$, so $L=H$.

Conversely, suppose $L=H$. Then $g\in N_G(H)$. By Theorem~\ref{thm:triple}, the product of two simple objects with stabilizer $H$ is again simple with stabilizer $H$, and their degrees multiply in $N_G(H)/H$. Let $k$ be the order of $gH$. The power $X^k$ then has stabilizer $H$ and lies in the component of the unit. This component is monoidally equivalent to $\mathbf{2Rep}(H)$ \cite[Lem.~3.2.5]{DY}. Under this equivalence, the simple objects with stabilizer $H$ correspond to actions of $H$ on $\Vect$. Their equivalence classes form the group $H^2(H;\A)$ under tensor product. Hence $X^k$ is invertible, and so is $X$.
\end{proof}

For $n\in N_G(H)$, the condition $[\xi_n]=[1]$ depends only on $nH$ by \eqref{eq:transport-cochain}. Set
$$
  \Gamma=\{nH\in N_G(H)/H:[\xi_n]=[1]\text{ in }H^3(H;\A)\}.
$$

\begin{proposition}\label{prop:invertible-extension}
Let $\operatorname{Inv}(\C)$ be the group of equivalence classes of invertible objects. Then $\Gamma$ is a subgroup of $N_G(H)/H$, and the degree map $[X(n,H,\phi)]\mapsto nH$ fits into an exact sequence
$$
  1\longrightarrow H^2(H;\A)\longrightarrow\operatorname{Inv}(\C)
    \xrightarrow{\deg} \Gamma\longrightarrow 1.
$$
\end{proposition}

\begin{proof}
By Lemma~\ref{lem:invertible-triples} and the fusion formula, taking degrees defines a group homomorphism $\operatorname{Inv}(\C)\to N_G(H)/H$. Its image is $\Gamma$, since $nH$ occurs precisely when $\dd\phi=\xi_n$ has a solution.

At degree $1$, the full-stabilizer labels are $X(1,H,\theta_1z)$ with $z\in Z^2(H;\A)$. Since $\ell_1^*\alpha=1$, changing the base point acts on $[z]$ by an inner automorphism of $H$, which acts trivially on $H^2(H;\A)$. Equation~\eqref{eq:Phi} gives multiplication by the product of cocycles, so $[z]\mapsto[X(1,H,\theta_1z)]$ identifies $H^2(H;\A)$ with the kernel.
\end{proof}

To compute the extension class, choose a representative $n_{\gamma}$ of each coset $\gamma$ and a lift $X_{\gamma}=X(n_{\gamma},H,\phi_{\gamma})$ for each $\gamma\in\Gamma$, with $n_1=1$ and $\phi_1=\theta_1$. For $\gamma,\eta\in\Gamma$, the product $X_{\gamma}\boxtimes_A X_{\eta}$ has a single summand, indexed by $x=1$, and the map $m$ in \eqref{eq:Phi} is $\operatorname{Ad}_{n_{\gamma}^{-1}}$. Let $\widehat\Phi_{\gamma,\eta}$ be the cochain obtained by applying \eqref{eq:Phi} and transporting the result from $n_{\gamma}n_{\eta}$ to $n_{\gamma\eta}$ by \eqref{eq:transport-cochain}.

The quotient $\widehat\Phi_{\gamma,\eta}/\phi_{\gamma\eta}$ is a 2-cocycle since the two cochains have the same differential. The factor set and the action of $\Gamma$ are
\begin{equation}\label{eq:extension-factor}
  c(\gamma,\eta)=[\widehat\Phi_{\gamma,\eta}/\phi_{\gamma\eta}]\in H^2(H;\A),
  \qquad
  \gamma\cdot[z]=[(\operatorname{Ad}_{n_{\gamma}^{-1}})^*z].
\end{equation}
Indeed, write $A_{[z]}=X(1,H,\theta_1z)$. The fusion formula gives $X_{\gamma}X_{\eta}=A_{c(\gamma,\eta)}X_{\gamma\eta}$ and $X_{\gamma}A_{[z]}=A_{\gamma\cdot[z]}X_{\gamma}$ on equivalence classes. Inner automorphisms of $H$ act trivially on $H^2(H;\A)$, so the action depends only on $\gamma$. Associativity gives the cocycle identity for $c$, and changing the lifts changes $c$ by a coboundary. Thus the extension splits if and only if $[c]$ vanishes in $H^2(\Gamma;H^2(H;\A))$.

\begin{corollary}\label{cor:graded-recognition}
The category $\C$ is monoidally equivalent to $\TV_{\mathcal{G}}^{\pi'}$ for some finite 2-group $\mathcal G$ and twist $[\pi']\in H^4(B\mathcal G;\A)$ if and only if $H$ is normal and abelian, $\Gamma=G/H$, and the extension in Proposition~\ref{prop:invertible-extension} splits.
\end{corollary}

\begin{proof}
By \cite[Rem.~2.1.17]{DR}, such a model exists precisely when every component contains an invertible object, the endomorphism category of the unit is pointed, and the degree map admits a group section. The components containing invertible objects are indexed by $\Gamma$. They exhaust $H\backslash G/H$ precisely when $H$ is normal and $\Gamma=G/H$. Since $\operatorname{End}_{\C}(\mathbf{1})\simeq\operatorname{Rep}(H)$ by \cite[Lem.~3.2.5]{DY}, the pointedness condition is equivalent to $H$ being abelian. The degree map admits a group section precisely when the extension in Proposition~\ref{prop:invertible-extension} splits.
\end{proof}

\subsection{Free bimodules and the example of \texorpdfstring{$S_3$}{S3}}
For $g\in G$, let $F_g=A\boxtimes(\Vect_g\boxtimes A)$ be the free bimodule generated in degree $g$. Writing $\bar g\in\mathcal R$ for the representative of $HgH$, we have $F_g\simeq X(\bar g,1,1)$. Both stabilizers in $F_gF_f$ are trivial, so each $x\in H$ gives a summand. Thus Theorem~\ref{thm:triple} recovers the formula $F_gF_f=\bigoplus_{x\in H}F_{gxf}$ of \cite[Prop.~3.3.1]{DY}.

Take $G=S_3$, $H=\langle s\rangle$ with $s=(12)$, and $\pi=\psi=1$. Write $r=(123)$. The double coset $H$ gives the unit and $X=F_1$. The other double coset $HrH$ gives a single simple object $D=F_r$, since $H_r=1$. In $D^2$, the terms indexed by $1$ and $s$ have degrees $r^2\in HrH$ and $rsr=s\in H$, respectively. The products are therefore
\begin{equation}\label{eq:s3}
 X^2=2X,\qquad XD=DX=2D,\qquad D^2=X+D .
\end{equation}
These are the surface fusion rules obtained by gauging the transposition subgroup of $S_3$ in \cite[Eq.~(119) and \S4.6]{BBST}.

\subsection{The nine fusion tables for \texorpdfstring{$C_2^3$}{C2 cubed}}\label{subsec:cubed}
Let $V=C_2^3$. We write $V=Q\oplus H$ additively. For $q\in Q\cong V/H$, let $\C_q$ be the component corresponding to the coset $g_q+H$, where $g_q=(q,0)$. In particular, $q=0$ labels the coset $H$ itself. In the tables, if the two inputs lie in $\C_q$ and $\C_r$, then every summand of their product lies in $\C_s$ with $s=q+r$. A line and a plane are subspaces of dimensions one and two over $\mathbb F_2$, respectively. We choose coordinates $x,y,z$ on $V$.

Table~\ref{tab:cubed-census} groups the 36 monoidal equivalence classes by their fusion rules. For each pair $(d,n)$ in the table, we choose the representative $\C(V,H_d,(-1)^{P_n},1)$. By the equivalence criterion of D\'ecoppet and Yu \cite[Prop.~3.4.4]{DY}, we may choose coordinates so that $H$ is the span $H_d$ of the last $d=\dim H$ coordinate vectors. The automorphisms of $V$ preserving $H_d$ then act on the classes $[\pi]\in H^4(V;\A)$ whose restriction to $H_d$ is trivial. In the basis of Appendix~\ref{app:cubed-calculations}, the eight binary coefficients of a quartic class give its integer label. In each orbit we take the least label $n$ and evaluate the corresponding homogeneous quartic $P_n\in\mathbb F_2[x,y,z]$ using ordered cup products. For each admissible label, every monomial of $P_n$ contains a coordinate that vanishes on $H_d$. Thus $\pi=(-1)^{P_n}$ restricts identically to $1$ on $H_d$, and any compatible $\psi$ is a three-cocycle. Pullback along $V\to H_d$ extends $\psi$ to a three-cocycle on $V$, so the same criterion allows us to take $\psi=1$. For fixed $d$ the monoidal equivalence classes correspond exactly to these orbits. Their numbers for $d=0,1,2,3$ are $7,19,9,1$, respectively.

When $H$ is a line, the subgroup $\Gamma\le Q$ of Proposition~\ref{prop:invertible-extension} consists of those $q$ for which $\C_q$ contains an invertible simple object. When $H$ is a plane, $Q\cong C_2$, so the two components are indexed by $q=0$ and $q=1$. The parameter $m$ counts the lines $L\le H$ that occur as stabilizers of simple objects in $\C_1$. If $m=3$, the coefficient $p\in\mathbb F_2$ in \eqref{eq:ranktwo-anomaly} determines the square of an invertible object in $\C_1$.

\begin{table}[htbp]
\centering\small
\caption{Fusion rules for $V=C_2^3$. The pair $(d,n)$ specifies a monoidal equivalence class. The rank is the number of equivalence classes of simple objects.}\label{tab:cubed-census}
\begin{tabular}{@{}cclcccl@{}}
\toprule
Case & $d$ & Parameters & Fusion rule & Rank & $\operatorname{Inv}(\C)$ & Quartic labels $n$\\
\midrule
1 & 0 & --- & $U_qU_r=U_s$ & 8  & $C_2^3$ & $0,1,6,8,10,44,99$\\
2 & 1 & $\Gamma=Q$       & \eqref{eq:cyclic-rule} & 8  & $C_2^2$ & $0,1,4,6,8,10,24,28,44$\\
3 & 1 & $[Q:\Gamma]=2$   & \eqref{eq:cyclic-rule} & 6  & $C_2$   & $64,65,66,80,81,82,96,97,98,99$\\
4 & 2 & $m=3,\ p=0$ & \eqref{eq:ranktwo-rule} & 12 & $C_2^2$ & $0,1$\\
5 & 2 & $m=3,\ p=1$ & \eqref{eq:ranktwo-rule} & 12 & $C_4$   & $8$\\
6 & 2 & $m=2$       & \eqref{eq:ranktwo-rule} & 9  & $C_2$   & $16,17$\\
7 & 2 & $m=1$       & \eqref{eq:ranktwo-rule} & 8  & $C_2$   & $2,6,10$\\
8 & 2 & $m=0$       & \eqref{eq:ranktwo-rule} & 7  & $C_2$   & $82$\\
9 & 3 & --- & \eqref{eq:cubed-burnside} & 30 & $C_2^3$ & $0$\\
\bottomrule
\end{tabular}
\end{table}

If the input stabilizers are $L,J\le H$, Theorem~\ref{thm:triple} gives one summand for each orbit of the diagonal $H$-action on $(H/L)\times(H/J)$. The difference of the two coset representatives, taken modulo $L+J$, identifies these orbits, so
\begin{equation}\label{eq:cubed-orbits}
 (H/L)\times(H/J)\simeq\coprod_{H/(L+J)}H/(L\cap J).
\end{equation}
Thus there are $[H:L+J]$ summands, each with stabilizer $L\cap J$. An integer in front of an object in a table denotes its multiplicity.

\textbf{Case 1.} When $H=0$, the unique simple object in $\C_q$ is $U_q=X(g_q,0,1)$. Its product with $U_r$ is $U_{q+r}$.

\textbf{Cases 2 and 3.}\label{subsec:13} Take $H=\langle e_z\rangle$ and $Q=\langle e_x,e_y\rangle$. We may write $P=h_4(x,y)+h_3(x,y)z+\ell(x,y)z^3$, where $h_j$ is homogeneous of degree $j$ and $\ell$ is linear. Equation~\eqref{eq:stabilizer-cocycle} gives $\xi_{g_q}=(-1)^{\ell(q)z^3}$. Each component $\C_q$ contains the free simple object $R_q=X(g_q,0,1)$. A simple object with stabilizer $H$ exists precisely when $\ell(q)=0$, since $(-1)^{z^3}$ represents the nontrivial class in $H^3(C_2;\A)$. In that case we take $I_q=X(g_q,H,1)$, and the triviality of $H^2(C_2;\A)$ makes it unique up to equivalence. Hence $\Gamma=\ker\ell$. The following formula gives case~2 when $\ell=0$ and case~3 when $\ell\ne0$.
\begin{equation}\label{eq:cyclic-rule}
\begin{array}{c|cc}
  & I_r & R_r\\ \hline
 I_q & I_s & R_s\\
 R_q & R_s & 2R_s
\end{array}
\qquad (R_q:\ q\in Q,\quad I_q:\ q\in \Gamma).
\end{equation}

\textbf{Cases 4--8.}\label{subsec:22} Now take $H=\langle e_y,e_z\rangle$ and identify $Q=\langle e_x\rangle$ with $\mathbb F_2$. The quartic representative has the form
\begin{equation}\label{eq:ranktwo-anomaly}
 P=x f_H+p x^2yz+x^3\lambda,\qquad
 f_H=ay^3+by^2z+cz^3,
\end{equation}
where $a,b,c,p\in\mathbb F_2$ and $\lambda\in H^*=\operatorname{Hom}_{\mathbb F_2}(H,\mathbb F_2)$. The cubic $f_H$ determines which stabilizers occur in $\C_1$. Indeed, $[\xi_{g_1}]=[(-1)^{f_H}]$, and on the three lines generated by $e_y,e_z,e_y+e_z$ its exponents are respectively $a,c,a+b+c$ times the cube of the nonzero character of the line. Thus $m$ counts the zeros among $a,c,a+b+c$. Any permutation of the three lines is induced by an automorphism of $H$, so their number determines the set of line labels up to relabeling.

For each $q\in\{0,1\}$, the component $\C_q$ contains the free simple object $R_q=X(g_q,0,1)$. For a line $L\le H$ we write $X_q^L=X(g_q,L,1)$. This object exists for every $L$ when $q=0$ and for the $m$ lines with zero exponent when $q=1$. In these cases $\xi_{g_q}|_L=1$ for our ordered cup representatives, and $H^2(L;\A)=0$. For stabilizer $H$ the two labels are
\[
 I_q^i=X(g_q,H,(-1)^{iyz}),\qquad i\in\mathbb F_2 .
\]
For $q=0$, both labels occur. For $q=1$, they occur precisely when $m=3$. Indeed, the three exponents vanish simultaneously exactly when $f_H=0$, and $(-1)^{yz}$ represents the generator of $H^2(H;\A)$. This gives six simple objects in $\C_0$. The component $\C_1$ contains $1+m$ simple objects with proper stabilizer, together with the two invertibles when $m=3$.

For each pair of labels that occurs, the product is given by
\begin{equation}\label{eq:ranktwo-rule}
\begin{array}{c|ccc}
 & I_r^j & X_r^J & R_r\\ \hline
 I_q^i & I_s^{\,i+j+pqr} & X_s^J & R_s\\
 X_q^L & X_s^L & \begin{cases}2X_s^L& L=J,\\ R_s& L\ne J\end{cases} & 2R_s\\
 R_q & R_s & 2R_s & 4R_s
\end{array}
\end{equation}
For a proper intersection of stabilizers, \eqref{eq:cubed-orbits} determines the product since the output cochain is unique up to equivalence. For full stabilizers, we use the middle representative $0$, so the base-point cochain is trivial. The comparisons \eqref{eq:theta-values} are also trivial in the components containing these objects. The only term of \eqref{eq:ranktwo-anomaly} with two quotient factors and two stabilizer factors is $p x^2yz$. Its contribution in \eqref{eq:betaexplicit} is $(-1)^{pqr yz}$ up to coboundary, which adds $pqr$ to the cocycle index in $\mathbb F_2$. In particular, when $m=3$, $(I_1^0)^2=I_0^p$. Here $I_0^0$ is the unit. This gives the groups $C_2^2$ and $C_4$ in cases~4 and~5, respectively.

\textbf{Case 9.} Here $H=V$, so $q=0$ is the only component label. We take $P=0$. By \cite[Lem.~3.2.5]{DY}, the whole fusion 2-category is then monoidally equivalent to $\mathbf{2Rep}(V)$. A simple object is given by a subspace $S\le V$ together with a class in $H^2(S;\A)$. Identify this group with $\bigwedge^2S^*$ by sending a cocycle $\phi$ to the alternating form $\omega$ determined by $\phi(u,v)/\phi(v,u)=(-1)^{\omega(u,v)}$. We write $\phi_\omega$ for any cocycle representing this form. For $S=0$ or a line $L$ there is only the trivial form, giving $R=X(0,0,1)$ and $X_L=X(0,L,1)$. For a plane $M$, let $\varepsilon_M$ be its unique nonzero alternating form. Then $Y_M^i=X(0,M,\phi_{i\varepsilon_M})$ for $i\in\mathbb F_2$. Finally, $U_\omega=X(0,V,\phi_\omega)$ for $\omega\in\bigwedge^2V^*$. These families have $1,7,14,8$ members, respectively.

The product adds the two forms after restriction to the intersection of the stabilizers \cite{GRR}. For example, two distinct planes meet in a line, where both forms restrict to zero. As their sum is $V$, \eqref{eq:cubed-orbits} gives a single copy of $X_{M\cap N}$. Writing $\eta|_M=\eta_M\varepsilon_M$, we obtain the following symmetric table, with only its upper half displayed.
\begin{equation}\label{eq:cubed-burnside}
\renewcommand{\arraystretch}{1.25}
\begin{array}{c|cccc}
 & R & X_J & Y_N^j & U_\eta\\ \hline
 R & 8R & 4R & 2R & R\\
 X_L & & \begin{cases}4X_L&L=J,\\2R&L\ne J\end{cases}
       & \begin{cases}2X_L& L\subset N,\\R&L\not\subset N\end{cases} & X_L\\
 Y_M^i & & & \begin{cases}2Y_M^{i+j}& M=N,\\ X_{M\cap N}&M\ne N\end{cases} & Y_M^{i+\eta_M}\\
 U_\omega & & & & U_{\omega+\eta}
\end{array}
\end{equation}
The rank and the isomorphism type of $\operatorname{Inv}(\C)$ distinguish all nine rows of Table~\ref{tab:cubed-census}.

\FloatBarrier

\subsection{Morita equivalence for \texorpdfstring{$C_2^3$}{C2 cubed}}\label{subsec:cubed-morita}
Every $\C(V,H,\pi,\psi)$ is Morita equivalent to $\TV_V^\pi$ through the module 2-category of $\Vect_H^\psi$-modules in $\TV_V^\pi$. By \cite[Prop.~3.4.1]{DY}, $\TV_V^\pi$ and $\TV_V^{\pi'}$ are Morita equivalent precisely when some automorphism of $V$ carries $[\pi']$ to $[\pi]$. Thus to pass from monoidal equivalence to Morita equivalence we forget the distinguished subgroup $H$. Starting from the quartic $P_n$ of a class $(d,n)$, we now take the least label under all of $\mathrm{GL}(V)$.

\begin{table}[htbp]
\centering\small
\caption{The seven Morita equivalence classes for $V=C_2^3$. Each entry lists the labels $n$ for the indicated $d$. Two monoidal equivalence classes are Morita equivalent precisely when their labels occur in the same row. The last column counts the monoidal equivalence classes in that row.}\label{tab:cubed-morita}
\begin{tabular}{@{}c l l c r@{}}
\toprule
$d=0$ & $d=1$ & $d=2$ & $d=3$ & Total\\
\midrule
$0$  & $0$             & $0$     & $0$ & 4\\
$1$  & $1,4,64$        & $1,2$   & --- & 6\\
$6$  & $6,65,96$       & $6$     & --- & 5\\
$8$  & $8,24,80$       & $8,16$  & --- & 6\\
$10$ & $10,28,66,81,98$& $10,17$ & --- & 8\\
$44$ & $44,82,97$      & $82$    & --- & 5\\
$99$ & $99$            & ---     & --- & 2\\
\bottomrule
\end{tabular}
\end{table}

Since every automorphism preserves the zero subgroup, the $d=0$ column gives the seven full orbit minima. For example, exchanging $x$ and $z$ interchanges $P_4=x^3z$ and $P_{64}=xz^3$ and moves the line $H=\langle e_z\rangle$ to $\langle e_x\rangle$. Thus $(1,4)$ and $(1,64)$ are Morita equivalent, although their fusion tables have ranks 8 and 6.

\FloatBarrier

\subsection{Finite 2-group models for \texorpdfstring{$C_2^3$}{C2 cubed}}\label{sec:explicit-models}
Let $\mathbb G(S,E,\nu)$ denote a skeletal 2-group with group of objects $S$, automorphism group $E$ of the unit, trivial $S$-action on $E$, and Postnikov class $\nu\in H^3(S;E)$. Thus $\nu=0$ gives the split 2-group $S[0]\times E[1]$. Write $\widehat S=\Hom(S,\A)$.

\begin{proposition}
Table~\ref{tab:cubed-models} gives a graded model for 19 classes and a representation model for 10 classes. The four classes $(d,0)$ admit both. For each row the models indicated are
\[
 \TV_{\mathbb G(Q,\widehat H,\nu_{\rm gr})}^{\,\mathrm{pr}^*[(-1)^{P_Q}]},
 \qquad
 \mathbf{2Rep}\bigl(\mathbb G(H,\widehat Q,\nu_{\rm rep})\bigr),
\]
respectively, where $\mathrm{pr}\colon B\mathbb G(Q,\widehat H,\nu_{\rm gr})\to BQ$ is the projection to the object group. The induced pullback $\mathrm{pr}^*:H^4(BQ;\A)\to H^4(B\mathbb G(Q,\widehat H,\nu_{\rm gr});\A)$ gives the twist of the graded model by evaluating $(-1)^{P_Q}$ on the four successive $Q$-labels.
\end{proposition}

The first $3-d$ coordinates of $V$ give the coordinates on $Q$, and the last $d$ give those on $H$. In the graded column of Table~\ref{tab:cubed-models}, let $T$ be a three-cocycle representing the indicated class in $H^3(Q;\mathbb F_2)$ and let $a\in H^*$. The term $T\otimes a$ in the Postnikov class is represented by the $\widehat H$-valued three-cocycle
\[
 (T\otimes a)(q_1,q_2,q_3)(h)
   =(-1)^{T(q_1,q_2,q_3)\,a(h)},
 \qquad q_1,q_2,q_3\in Q,\quad h\in H.
\]
Thus a triple on which $T$ takes the value $1$ gives the character $h\mapsto(-1)^{a(h)}$, and a triple on which $T$ vanishes gives the trivial character. In the representation column the same formula applies with $H$ and $Q$ interchanged. For a graded model, $P_Q$ specifies the four-cocycle twist, with $P_Q=0$ giving the untwisted model.

\begin{table}[!htbp]
\centering\small
\caption{Finite 2-group models for 25 classes $(d,n)$. The object and unit automorphism groups are $Q,\widehat H$ for a graded model and $H,\widehat Q$ for a representation model. A dash means that the class admits no model of that type.}\label{tab:cubed-models}
\begin{tabular}{@{}rr ll l@{}}
\toprule
&&\multicolumn{2}{c}{Graded model}&\multicolumn{1}{c}{Representation model}\\
\cmidrule(lr){3-4}\cmidrule(lr){5-5}
$d$ & $n$ & $\nu_{\rm gr}$ & $P_Q$ & $\nu_{\rm rep}$\\
\midrule
0&0&$0$&$0$&$0$\\
0&1&$0$&$P_1$&---\\
0&6&$0$&$P_6$&---\\
0&8&$0$&$P_8$&---\\
0&10&$0$&$P_{10}$&---\\
0&44&$0$&$P_{44}$&---\\
0&99&$0$&$P_{99}$&---\\
\midrule
1&0&$0$&$0$&$0$\\
1&1&$0$&$x^3y$&---\\
1&4&$x^3\otimes z$&$0$&---\\
1&6&$x^3\otimes z$&$xy^3$&---\\
1&8&$x^2y\otimes z$&$0$&---\\
1&10&$x^2y\otimes z$&$xy^3$&---\\
1&24&$(x^2y+xy^2)\otimes z$&$0$&---\\
1&28&$(x^3+x^2y+xy^2)\otimes z$&$0$&---\\
1&44&$(x^3+x^2y+y^3)\otimes z$&$0$&---\\
1&64&---&---&$z^3\otimes x$\\
\midrule
2&0&$0$&$0$&$0$\\
2&1&$x^3\otimes y$&$0$&---\\
2&2&---&---&$y^3\otimes x$\\
2&8&---&---&$yz(y+z)\otimes x$\\
2&10&---&---&$y(y^2+yz+z^2)\otimes x$\\
2&16&---&---&$y^2z\otimes x$\\
2&82&---&---&$(y^3+y^2z+z^3)\otimes x$\\
\midrule
3&0&$0$&$0$&$0$\\
\bottomrule
\end{tabular}
\end{table}

\begin{proof}
All subgroups of $V$ are normal and abelian. By Corollary~\ref{cor:graded-recognition}, a graded model exists precisely when every component contains an invertible object and the degree map $\operatorname{Inv}(\C)\to Q$ splits. In Table~\ref{tab:cubed-census} these are cases 1, 2, 4 and 9. Case 5 has invertibles for both $q=0$ and $q=1$, but its degree map $C_4\to C_2$ has no section.

For the 19 classes in these four cases, write $P_n=P_Q+\sum_i a_iT_i$ in mod-two cohomology, where $P_Q$ depends only on $Q$, each $a_i$ is linear on $H$, and each $T_i$ is cubic on $Q$. Ordering the mixed cups as $a_i\smile T_i$ gives the character pairing of D\'ecoppet and Yu \cite[Eq.~(12), Prop.~4.2.1]{DY}, with Postnikov class $\nu_{\rm gr}=\sum_i T_i\otimes a_i$. Any three-cochain $c$ comparing the two orders restricts to a cocycle on $H$, since both representatives restrict to $1$ there. For the projection $r:V\to H$, the cochain $c/r^*(c|_H)$ has the same coboundary and restricts to $1$, preserving $\psi=1$.

Inflating $(-1)^{P_Q}$ along $V\to Q$ contributes the scalar $1$ to every bimodule or balancing constraint with an $H$-input. On four homogeneous bimodules of quotient degrees $q_1,\ldots,q_4$, it multiplies the pentagonator by $(-1)^{P_Q(q_1,\ldots,q_4)}$. Under the graded equivalence this gives precisely the twist $\mathrm{pr}^*[(-1)^{P_Q}]$. For example, $(1,6)$ has $P_6=x^3z+xy^3$. Its graded model has object group $Q=C_2^2$ and unit automorphism group $\widehat H=C_2$. The term $x^3z$ gives the Postnikov class $x^3$ with values in this $C_2$, and the twist is the pullback of $[(-1)^{xy^3}]$ from $BQ$.

For a representation model, the same pairing interchanges the roles of $H$ and $Q$ \cite[Cor.~4.2.3]{DY}. Indeed, a 2-group with object group $S$ and unit automorphism group $E$ gives a presentation with ambient group $S\ltimes\widehat E$ and distinguished subgroup $S$. A monoidal equivalence must identify this pair with $(V,H)$ \cite[Prop.~3.4.4]{DY}. It follows that $S\cong H$, $E\cong\widehat Q$, and the action is trivial. Thus a representation model exists exactly when the orbit of $[(-1)^{P_n}]$ under the automorphisms preserving $H$ contains a class $[(-1)^P]$ with $P$ a sum of terms $aT$, where $a$ is linear on $Q$ and $T$ is cubic on $H$. The sum $\sum T\otimes a$ then gives $\nu_{\rm rep}$.

For $H=0$ or $V$ the only such quartic is zero. When $H$ is a line, the nonzero possibilities have the single representative $xz^3$, giving $(1,64)$. When $H$ is a plane, the quartic has the form $xT(y,z)$. The nonzero binary cubics fall into five orbits under change of basis on $H$, represented by the entries with $\nu_{\rm rep}\ne0$ for $d=2$. They correspond, respectively, to a cube, three distinct linear factors, a linear factor times an irreducible quadratic, a double and a simple linear factor, and an irreducible cubic, in the order $n=2,8,10,16,82$. For instance, at $(2,8)$ the representative $P_8=x^2yz$ and the quartic $x\,yz(y+z)$ differ by $\mathrm{Sq}^1(xyz)$, so they define the same sign class. This row therefore gives the representations of the 2-group with object group $C_2^2$, unit automorphism group $C_2$, and Postnikov class $yz(y+z)$.
\end{proof}

The eleven classes admitting neither model are $(2,6)$, $(2,17)$, and the nine line classes
$$
 (1,n),\qquad n\in\{65,66,80,81,82,96,97,98,99\}.
$$
Among the ten classes with the fusion rules of case 3, only \((1,64)\) admits a finite 2-group model, namely the representation model in Table~\ref{tab:cubed-models}.

\FloatBarrier

\appendix
\section{Comparison with the stabilizer cocycle}\label{app:comparison}
Let $s$ be a simple object of degree $g$ in the category underlying a bimodule $M$, and let $\Delta_gL$ be the stabilizer of its isomorphism class. For $a_i\in L$, the pair $(a_i,b_i)$ with $b_i=g^{-1}a_i g$ fixes the class of $s$. Lemma~\ref{lem:comparison-scalar} uses $q_i:a_i s\to s b_i$. In the stabilizer convention of \cite[Prop.~3.2.1]{DY}, the right action is by $b_i^{-1}$. Applying it to $q_i$ and then using the right module associator $r^M$ gives
\[
 (a_i s)b_i^{-1}\xrightarrow{\,q_i\,b_i^{-1}\,}(s b_i)b_i^{-1}\xrightarrow{\,r^M\,}s(b_i b_i^{-1})\cong s .
\]

For two such pairs the composite in Lemma~\ref{lem:comparison-scalar} applies $q_1$ and then $q_2$, leaving the right label $b_2b_1$. If we pass to stabilizer maps at each step, we cancel $b_1 b_1^{-1}$ first and $b_2 b_2^{-1}$ second. If we compose the maps of the lemma first, we cancel $(b_2b_1)(b_2b_1)^{-1}$ in one step. Since $(b_2b_1)^{-1}=b_1^{-1}b_2^{-1}$, the successive cancellations on the right give the composite
\[
 (b_2b_1)(b_1^{-1}b_2^{-1})\longrightarrow b_2\bigl(b_1(b_1^{-1}b_2^{-1})\bigr)\longrightarrow b_2\bigl((b_1b_1^{-1})b_2^{-1}\bigr)\longrightarrow b_2b_2^{-1}\longrightarrow 1 .
\]
With the associator convention of Section~\ref{sec:classification}, the first arrow contributes $\psi(b_2,b_1,(b_2b_1)^{-1})^{-1}$ and the second contributes $\psi(b_1,b_1^{-1},b_2^{-1})$. Write $D(b_2,b_1)$ for their product. Reassociating past the object of degree $g$ also contributes factors of $\pi$.

The two products can be compared on the whole group $H_g$ before we restrict to $L$. $(\ell_g^*\alpha)^{-1}$ is the associator in the convention of Lemma~\ref{lem:comparison-scalar}, where the two-step composite has coefficient $b_S^{-1}$ relative to the direct map. To express this product using the stabilizer associator $\xi_g$, we multiply its coefficient by $\theta_g(a_2,a_1)$. For three successive pairs $(a_i,b_i)$, compatibility with the associators requires
\begin{equation}\label{eq:theta-coherence}
 (\ell_g^*\alpha)^{-1}(a_3,a_2,a_1)\,\theta_g(a_2,a_1)\theta_g(a_3,a_2a_1)=\xi_g(a_3,a_2,a_1)\,\theta_g(a_3,a_2)\theta_g(a_3a_2,a_1).
\end{equation}
For $a_1,a_2\in H_g$ we take the normalized two-cochain
\begin{equation}\label{eq:theta-values}
\begin{aligned}
\theta_g(a_2,a_1)
 &=\frac{\pi(a_2,a_1,a_1^{-1}g,b_2^{-1})\,\pi(a_2a_1,g,b_1^{-1},b_2^{-1})}
         {\pi(a_2g,b_1,b_1^{-1},b_2^{-1})\,\pi(a_2,g,b_1,b_1^{-1})}\\
 &\qquad{}\cdot\pi(g,b_2,b_1,(b_2b_1)^{-1})\,\pi(a_2,a_1,g,b_1^{-1})\,
          \frac{\psi(b_1,b_1^{-1},b_2^{-1})}{\psi(b_2,b_1,(b_2b_1)^{-1})}.
\end{aligned}
\end{equation}

To verify \eqref{eq:theta-coherence}, write $M_g(a_2,a_1)$ for the product of the six $\pi$-factors in \eqref{eq:theta-values}, so that $\theta_g(a_2,a_1)=M_g(a_2,a_1)D(b_2,b_1)$. Applying $\dd\psi=\pi|_H$ to the three quadruples in the following expression for $R$ gives
\[
\begin{aligned}
 R&=\frac{\pi(b_1,b_1^{-1},b_2^{-1},b_3^{-1})}{\pi(b_3,b_2,b_1,(b_3b_2b_1)^{-1})\,\pi(b_2,b_1,(b_2b_1)^{-1},b_3^{-1})},\\
 (\dd D)(b_3,b_2,b_1)&=R\,\frac{\psi(b_3,b_2,b_1)}{\psi(b_1^{-1},b_2^{-1},b_3^{-1})}.
\end{aligned}
\]
Substituting the six factors of $M_g$ into $(\dd M_g)(a_3,a_2,a_1)$ and applying $\dd\pi=1$ gives
\[
\begin{aligned}
 R\,(\dd M_g)(a_3,a_2,a_1)
 &=\frac{\pi(a_1,a_1^{-1}g,b_2^{-1},b_3^{-1})\,\pi(g,b_1^{-1},b_2^{-1},b_3^{-1})}
         {\pi(a_2,a_1,(a_2a_1)^{-1}g,b_3^{-1})\,\pi(a_3,a_2,a_1,(a_3a_2a_1)^{-1}g)}\\
 &\qquad{}\cdot\frac{\pi(a_3,a_2,a_1,g)\,\pi(a_3,g,b_2,b_1)}{\pi(a_3,a_2,g,b_1)\,\pi(g,b_3,b_2,b_1)} .
\end{aligned}
\]
The first fraction contains the $\pi$-factors of $\xi_g$ in \eqref{eq:stabilizer-cocycle}, and the second contains those of $\ell_g^*\alpha$ in \eqref{eq:alphaexpanded}. Multiplication by $(\dd D)(b_3,b_2,b_1)/R$ supplies the ratio of the right associators. Hence $\dd\theta_g=\xi_g\,\ell_g^*\alpha$, as required. Restricting to $L$ now gives the multiplication coefficient $\theta_g|_L/b_S|_{s,L}$ used in Theorem~\ref{thm:decorated-classification}.

\section{Quartic representatives for the two orbit calculations}\label{app:cubed-calculations}
We describe how to recover the quartic $P_n$ from a label in Table~\ref{tab:cubed-census} and how to find the label of a given presentation. Write $e_x,e_y,e_z$ for the chosen basis of $V$ and $x,y,z$ for the dual basis of $V^*$. Then $H^\bullet(V;\mathbb F_2)=\mathbb F_2[x,y,z]$. Every class in $H^4(V;\A)$ has a representative $(-1)^P$ for a homogeneous quartic $P$. Indeed, the coefficient sequence
\[
  1\longrightarrow\{\pm1\}\longrightarrow\A
    \xrightarrow{w\mapsto w^2}\A\longrightarrow1
\]
identifies the image of the sign map with the kernel of squaring on $H^4(V;\A)$. By the K\"unneth formula, $H^4(V;\A)\cong C_2^8$, so this kernel is the whole group.

For a homogeneous cubic $h$, the sign of $\mathrm{Sq}^1h$ is a coboundary. In fact, if $\widetilde h$ is an integral lift of a cocycle representing $h$, then $\dd(\mathrm i^{\widetilde h})=(-1)^{\dd\widetilde h/2}$, where $\dd\widetilde h/2\bmod 2$ represents $\mathrm{Sq}^1h$. Applying $\mathrm{Sq}^1$ to the cubic monomials gives seven independent relations among the quartics. In the image of the sign map, the fourth powers and the products of two squares vanish, and
\[
  xyz^2=x^2yz+xy^2z .
\]
These relations span the kernel of the sign map, which has dimension $15-8=7$. The signs of the following quartics give an ordered basis of $H^4(V;\A)$ over $\mathbb F_2$.
\[
 (x^3y,\ xy^3,\ x^3z,\ x^2yz,\ xy^2z,\ y^3z,\ xz^3,\ yz^3).
\]
If a class has coefficients $b_0,\ldots,b_7$ in this basis, write $n=\sum_{j=0}^7 2^j b_j$ and denote the representative by $P_n$. For example, $P_{10}=xy^3+x^2yz$.

Let $H$ be the span of the last $d$ coordinate vectors. First retain the labels for which $P_n$ restricts to zero on $H$ in the quartic quotient. For $d=0,1$ all 256 labels qualify. For $d=2$ the coefficients of $y^3z$ and $yz^3$ must vanish, leaving 64 labels. For $d=3$ only $n=0$ remains. To find the orbit of a retained label, substitute the row forms of each matrix in $\mathrm{GL}(V)_H=\{A\in\mathrm{GL}(V)\mid A(H)=H\}$ for $x,y,z$ in $P_n$, and reduce by the seven relations. Record the least label in the orbit and remove that orbit from the list. Repeating this calculation gives the $7,19,9,1$ representatives for $d=0,1,2,3$ listed in Table~\ref{tab:cubed-census}. For instance, when $d=2$, exchanging $y$ and $z$ takes $P_4=x^3z$ to $P_1=x^3y$, so both presentations represent the monoidal equivalence class labeled by $(2,1)$.

To determine the Morita equivalence class of a given pair $(d,n)$, we forget the subgroup $H_d$ and take the orbit of $[(-1)^{P_n}]$ under the full group $\mathrm{GL}(V)$. For each of the 168 elements of $\mathrm{GL}(V)$, we make the corresponding linear substitution in $P_n$ and reduce the result using the relations above. If $n'$ is the smallest label obtained, then the category labeled by $(d,n)$ is Morita equivalent to the one labeled by $(0,n')$. Thus $(d,n)$ belongs to the row of Table~\ref{tab:cubed-morita} with $n'$ in the $d=0$ column.

\end{document}